\documentclass[12pt]{article}

\usepackage{amssymb,amsmath,amsthm,graphics}
\newtheorem{defi}{Definition}
\newtheorem{thr}{Theorem}
\newtheorem{exa}[defi]{Example}
\newtheorem{lem}[defi]{Lemma}
\newtheorem{cor}[defi]{Corollary}
\newtheorem{pro}[defi]{Proposition}

\newcommand{\eps}{\varepsilon}
\renewcommand{\phi}{\varphi}
\newcommand{\A}{AP}
\newcommand{\F}{\mathcal{F}}
\renewcommand{\H}{\mathcal{H}}
\newcommand{\U}{\mathfrak{U}}
\renewcommand{\H}{\mathcal{H}}
\newcommand{\T}{\mathbb{T}}
\newcommand{\C}{\mathbb{C}}
\newcommand{\Z}{\mathbb{Z}}
\newcommand{\N}{\mathbb{N}}
\newcommand{\R}{\mathbb{R}}
\renewcommand{\b}{{}^{\flat}\!}
\newcommand{\Eins}{\mbox{1\hspace{-.9mm}\mbox{I}}}

\begin{document}
\title{Filters and subgroups associated \\
       with Hartman measurable functions}
\author{Gabriel Maresch\footnote{The author would like to thank the FWF for
        financial support through grant S8312}}
\maketitle
\begin{abstract}
A bounded function $\phi: G\to \C$ on an LCA group $G$ is called Hartman measurable 
if it can be extended to a Riemann integrable function $\phi^*: X\to \C$ on some group
compactification $(\iota_X,X)$, i.e. on
a compact group $X$ such  that $\iota_X: G\to X$ is 
a continuous homomorphism with image $\iota_X(G)$ dense in $X$ and
$\phi=\phi^*\circ\iota_X$. The concept of Hartman measurability of functions is a generalization
of Hartman measurability of sets, which was introduced - with different
nomenclature -  by S. Hartman
to treat number theoretic problems arising 
in diophantine approximation and equidistribution.
We transfer certain results concerning Hartman sets to this more
general setting. In particular we assign to each Hartman measurable function
$\phi$ a filter $\F(\phi)$ on $G$ and a subgroup $\Gamma(\phi)$ of the dual $\hat{G}$ and
show how these objects encode information about the involved
group compactification. We present
methods how this information can be recovered.\par\vspace{1em}

\par\vspace{1em}
\noindent 2000 Mathematics Subject Classification
37A45, %Relations with number theory and harmonic analysis 
43A60, %Almost periodic functions on groups and semigroups and their generalizations
(11K70)%Harmonic analysis and almost periodicity
\end{abstract}

\newpage \section{Introduction}
\subsection{Motivation:} In \cite{FP} the investigation of finitely additive measures in number 
theoretic context led to the concept of a Hartman measurable subset $H\subseteq G$ of
a discrete abelian group $G$. By definition, $H$ is Hartman
measurable if it is the preimage $H=\iota_X^{-1}(M)$ of a continuity set $M\subseteq X$
in  a group compactification $(\iota_X,X)$ of $G$. This, more explicitly, means that $\iota_X: G\to C$
is a group homomorphism with $\iota(G)$ dense in the compact group $X$ and that
$\mu_X(\partial M)=0$. Here $\partial M$ denotes the topological boundary of $M$ and
$\mu_X$ the normalized Haar measure on $X$. By putting $m_G(H)=\mu_X(M)$ one can define a 
finitely additive measure on the Boolean set algebra of all Hartman measurable
sets in $G$. For the special case $G=\Z$ a Hartman set $H\subseteq \Z$, by identification
with its characteristic function, can be considered to be a two-sided infinite
binary sequence, called a Hartman sequence. Certain number theoretic, ergodic and 
combinatorial aspects of Hartman sequences have been studied in \cite{SchW} and \cite{StW},
while \cite{Wi} presents a method to reconstruct the group compactification
$(\iota_X,X)$ for given $H$.

In order to benefit from powerful tools from functional and harmonic analysis it
is desirable to study appropriate generalizations of Hartman measurable sets by replacing
their characteristic functions by complex valued functions not only taking
the values 0 and 1. The natural definition of a Hartman measurable function $\phi: G \to \C$
is the requirement $\phi^*\circ\iota_X$, where $(\iota_X,X)$ is a group compactification
of $G$ and $\phi^*$ is integrable in the Riemann sense, i.e. its points of
discontinuity form a null set with respect to the Haar measure on $X$. This 
definition is equivalent to the one of R-almost periodicity, introduced
in \cite{H} by S. Hartman.  The investigation
of the space $\H(G)$ of all Hartman measurable functions on $G$ is the content of \cite{MW}.
Here we are going to
transfer ideas from \cite{Wi} into this context. Thus our
main topic is to describe $(\iota_X,X)$ only in terms of $\phi$. In particular we establish
further connections to Fourier analysis. The natural framework for our investigation
is that of LCA (locally compact abelian) groups. 

\subsection{Content of the paper} After the introduction we collect in section 2 the necessary
preliminary definitions and facts about Hartman measurable sets and functions.

Section 3 treats the following situation: Given a Hartman measurable function
$\phi: G\to X$ on an LCA group $G$, we know by the very definition of Hartman 
measurability that there is some group compactification $(\iota_X,X)$ of $G$ such that
$\phi=\phi^*\circ\iota_X$ for some Riemann integrable function $\phi^*: X\to \C$.
We say that $\phi$ can be realized in $(\iota_X,X)$ resp. by $\phi^*$.
It is easy to see
that in this case $\phi$ can be realized as well on any "bigger" compactification
$(\iota_{\tilde{X}},\tilde{X})$. The notion of "bigger" and "smaller"
is made more precise in the next section.

In particular every Hartman measurable function can be realized 
in the maximal group compactification of $G$, the Bohr compactification $(\iota_b,bG)$.
The question arises if there is
a realization of $\phi$ in a group compactification that is as "small" as possible. 
If a Hartman measurable function $\phi$ possesses a so called aperiodic realization
then the group compactification on which this aperiodic realization can be obtained 
is minimal (Theorem 1).
This approach works for arbitrary Hartman measurable $\phi$
if one allows "almost realizations", i.e. if one demands $\phi=\phi^*\circ\iota_X$ 
almost everywhere with respect to the finitely additive Hartman measure $m_G$ on $G$ rather
than $\phi=\phi^*\circ\iota_X$ everywhere on $G$ (Theorem \ref{weilthr}). Whenever $\phi$ is
even almost periodic one can guarantee $\phi=\phi^*\circ\iota_X$ everywhere on $G$.
The group compactification on which the minimal realization of $\phi$ occurs is
unique up to
equivalence of group compactifications. 
It can be obtained by a method involving filters on $G$ similar to that
presented in \cite{Wi}.

The content of section 4 is motivated by the following reasoning: Every group compactification
of the LCA group $G$ corresponds to a (discrete) subgroup $\Gamma$ of the dual $\hat{G}$ in
such a way that it is equivalent to the group compactification $(\iota_{\Gamma}, C_{\Gamma})$ defined by
$\iota_{\Gamma}: g\to (\chi(g))_{\chi \in \Gamma}$, $C_{\Gamma}:=
\overline{\iota_{\Gamma}(G)} \le \T^{\Gamma}$.
If $(\iota_X, X)$ is a group
compactification admitting an aperiodic almost realization of the Hartman measurable
function $\phi$, the corresponding subgroup $\Gamma\le\hat{G}$ contains all characters $\chi$ 
such that the corresponding Fourier coefficient $m_G(\phi\cdot\overline{\chi})$ does not
vanish. If $\phi$ is almost periodic or if $\phi$ can be realized on a finite dimensional
compactification this result is sharp in the sense that the subgroup ${\Gamma}$ is minimal
with the above property (Theorem \ref{mincorap}). For general Hartman measurable
functions the situation is more difficult. This is discussed and illustrated by an example.

Section 5 summarizes the main results and includes an illustrating diagram.

\section{Preliminaries and Notation} Throughout this paper
$G$ denotes always an LCA (locally compact abelian) group.
For group compactifications of $G$ let us 
write $(\iota_{X_1}, X_1)\le (\iota_{X_2}, X_2)$
iff there is a continuous group homomorphism $\pi: X_2\to X_1$ such that the diagram
\begin{center}
\includegraphics*{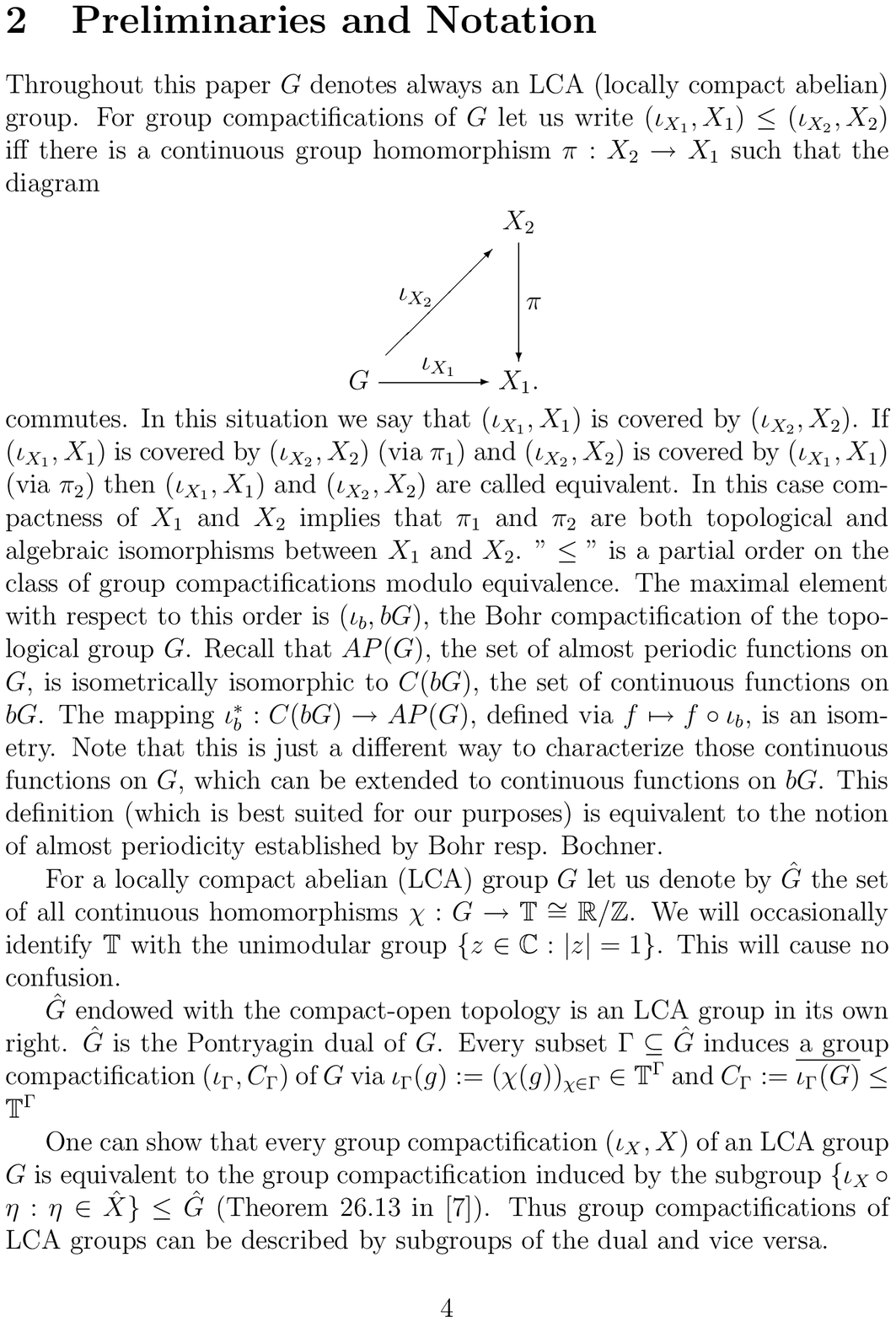}
\end{center}
commutes. In this situation we say that 
$(\iota_{X_1}, X_1)$ is covered by $(\iota_{X_2}, X_2)$.
If $(\iota_{X_1}, X_1)$ is covered by $(\iota_{X_2}, X_2)$ (via $\pi_1$)
and $(\iota_{X_2}, X_2)$ is covered by $(\iota_{X_1}, X_1)$ (via $\pi_2$)
then $(\iota_{X_1}, X_1)$ and $(\iota_{X_2}, X_2)$ are called 
equivalent. In this case compactness of $X_1$ and $X_2$ implies
that $\pi_1$ and $\pi_2$
are both topological and algebraic isomorphisms between $X_1$ and $X_2$.
$"\le"$ is a partial order on
the class of group compactifications 
modulo equivalence. The maximal element with respect to this order is $(\iota_b,bG)$,
the Bohr compactification of the topological
group $G$. Recall that $AP(G)$, the set of almost periodic functions on $G$,
is isometrically isomorphic to $C(bG)$, the set of continuous functions
on $bG$. The mapping $\iota_b^*: C(bG)\to AP(G)$, defined via
$f\mapsto f\circ\iota_b$, is an isometry. Note that this is just a different way to characterize those continuous functions
on $G$, which can be extended to continuous functions on $bG$. 
This definition (which is best suited for our purposes) is equivalent
to the notion of almost periodicity established by Bohr resp. Bochner.

For a locally compact abelian (LCA) group $G$ let us denote by $\hat{G}$ the set
of all continuous homomorphisms $\chi: G \to \T\cong \R/\Z$. We will occasionally identify
$\T$ with the unimodular group $\{z\in \C: |z|=1\}$. This will cause no confusion.

$\hat{G}$ endowed
with the compact-open topology is an LCA group in its own right. $\hat{G}$
is the Pontryagin dual of $G$. Every subset $\Gamma\subseteq \hat{G}$ induces a
group compactification $(\iota_{\Gamma},C_{\Gamma})$ of $G$ via
$\iota_{\Gamma}(g):=(\chi(g))_{\chi\in \Gamma}\in 
\T^{\Gamma}$ and
$C_{\Gamma}:=\overline{\iota_{\Gamma}(G)} \le \T^{\Gamma}$

One can show that every group compactification $(\iota_X,X)$ of an LCA group $G$
is equivalent to the group compactification induced by the subgroup
$\{\iota_X\circ\eta: \eta\in\hat{X}\}\le \hat{G}$ (Theorem 26.13 in \cite{HR}).
Thus group compactifications of LCA groups
can be described by subgroups of the dual and vice versa.

The system $\Sigma(G)\subseteq \mathfrak{P}(G)$ of all 
Hartman measurable sets on $G$, i.e. the system
of all preimages $\iota_b^{-1}(M)$ of $\mu_b$-continuity
sets in the Bohr compactification $(\iota_b,bG)$ of $G$,
is a Boolean set algebra and enjoys the
property that there exists a unique translation invariant finitely additive
probability measure $m_G$ on $\Sigma(G)$: 
$m_G(\iota_X^{-1}(M^*))=\mu_X(M^*)$. For details we refer to \cite{FP}.

Let us denote by $\Delta$ the symmetric difference of sets and by $\tau_g$
the translation operator on an abelian group defined by $\tau_g (h):= g+h$.
We introduce two mappings:
\begin{itemize}
\item for a Hartman measurable set $M$
denote by $d_M: G \to [0,1]$ the mapping $g\mapsto m_G( M \Delta \tau_g M )$,
\item for a $\mu_X$-continuity set $M^*$ on some group compactification $(\iota_X,X)$
denote by $d_{M^*}: X \to [0,1]$ the mapping $g\mapsto \mu_X( M^* \Delta \tau_g M^* )$.
\end{itemize}

Note that the mapping $d_{M^{*}}$ (and similarly the mapping
$d_{M}$)
can be used to define a translation invariant pseudometric by letting
$\rho_{M^{*}}(g,h):=d_{M^{*}}(g-h)$. The set of zeros $\{g: d_{M^{*}}(g)=0\}$ is always a
closed subgroup. We will denote this subgroup by $\ker d_{M^{*}}$.

Now consider sets of the form
$F(M,\eps):=\{g\in G: d_M(g)<\eps\}$
and denote by $\F(M)$ the filter on $G$ generated by $\{F(M,\eps):\eps > 0\}$, i.e.
the set of all $F\subseteq G$ such that there exists an $\eps>0$ with $F(M,\eps)\subseteq F$.
When we have a realization $M^*$ of $M$ on some group compactification $(\iota_X, X)$
we can transfer the topological data encoded in the neighborhood filter of the unit $0_X$ in $X$
to $G$ by considering the pullback induced by $\iota_X$. 

To be precise: Let $(\iota_X,X)$ be a group compactification and $\U(X,0_X)$
the filter of all neighborhoods
of the unit $0_X$ in $X$. By $\U_{(\iota_X,X)}$ we denote the filter on $G$
generated by $\iota_X^{-1}\left(\U(X,0_X)\right)$.
Note that if the mapping $\iota_X$ is one-one, $\iota_X^{-1}\left(\U(X,0_X)\right)$ is already a filter.

For $\Z$, the group of integers,
Theorem 2 in \cite{Wi} states that for any Hartman set $M\subseteq \Z$ 
there is a group compactification
$(\iota_X,X)$ such that $\F(M)$ and $\U(X,0_X)$ coincide. Hence the 
filter $\F(M)$ on $\Z$ contains much information about
the group compactification $(\iota_X,X)$:
{\em If $M\subseteq \Z$ is a Hartman measurable set and $(\iota_X, X)$ is a group compactification of
the integers such that $M$ can be realized on $X$ via the continuity set $M^*$ 
then $H= \ker d_{M^*}$ is a closed subgroup of $X$ and
$\F(M)=\U_{(\pi_H\circ \iota_X,X/H)}$, for $\pi_H: X\to X/H$ the canonical quotient mapping.}

In what follows we need to generalize this result to arbitrary (LCA) groups. This 
poses no problem since the proof 
given in \cite{Wi} for $G=\Z$ applies
verbatim to an arbitrary topological group.

Recall that for a filter $\F\subseteq \mathfrak{P}(X)$ on some set $X$
and a function $\phi: X\to \C$ the filter-limit 
$\F\!\!-\!\!\lim_{x\in X} \phi(x)$ is defined to be the unique $\lambda\in\C$ such
that $\forall \eps> 0$ we have $\{x\in X: |\phi(x)-\lambda|<\eps\}\in \F$.
In \cite{Wi} the filter $\F=\F(M)$ is also used to define the subgroup Sub$(M)$ of
$\T$ consisting of all those elements $\alpha$ such that
$\F\!\!-\!\!\lim_{n\in\Z} \lfloor n\alpha\rfloor=0$ 
(or equivalently: $\F\!\!-\!\!\lim_{n\in\Z} e^{2\pi i n\alpha}=1$).

All three objects - filter, compactification and subgroup - carry the same
information regarding a fixed Hartman set $M$. It is interesting to note
that any subgroup of a compact abelian group $G$ can be written as
$\{g\in G: \F\!\!-\!\!\lim_{\chi\in\hat{G}} \chi(g)=1\}$ for some filter
$\F$ on $\hat{G}$ (cf. \cite{BS}).
 
We transfer these concepts into our more general context.
To that cause we need the following definitions. Recall that 
a bounded function $f$ on a group compactification $(\iota_X,X)$
is called Riemann integrable iff the set $\mbox{disc}(f)$ of points of discontinuity
is a $\mu_X$-null set, for $\mu_X$ the normalized Haar measure on $X$.
Let us denote the set of all such functions by $R_{\mu_X}(X)$ or,
simply, $R(X)$. We use the
following characterization, a proof of which can be found in \cite{Ta}.

\begin{pro}\label{riemannchar} 
Let $X$ be a compact space and $\mu_X$ a finite positive Borel measure
on $X$. For a bounded real-valued $\mu_X$-measurable function $f$ the following
assertions are equivalent:
\begin{enumerate}
\item $f$ is Riemann integrable.
\item For every $\eps>0$ there exist continuous functions $g_{\eps}$ and 
$h_{\eps}$ such that $g_{\eps}\le f \le h_{\eps}$ and
$\int_X (h_{\eps}-g_{\eps}) d\mu_X< \eps$.
\end{enumerate}
\end{pro}

Let $\phi$ be a function defined on a topological group $G$ and $(\iota_b,bG)$ the
Bohr compactification of $G$. We call 
a function ${\phi}_b$ defined on $bG$
an extension (or realization) of $\phi$ iff $\phi={\phi}_b\circ\iota_b$.
For example: The set of almost periodic functions on $G$ coincides with 
the set of those functions that can be extended to continuous functions on the Bohr
compactification.

\begin{defi} Let $(\iota_b,bG)$ be the Bohr compactification of the topological group $G$.
We call a bounded function $\phi$ on $G$ Hartman measurable
iff $\phi$ can be extended to a Riemann integrable function $\phi_b$ on $bG$.
The set of Hartman measurable functions $\{\phi^*\circ\iota_b: \phi^*\ \in R_{\mu_b}(bG)\}$
is denoted by $\H(G)$.
\end{defi}

Given a Hartman measurable
function $\phi$, we say that \emph{$\phi^*$ realizes $\phi$} if $\phi^*$ is a Riemann integrable
function defined on some
group compactification $(\iota_X, X)$ such that $\phi=\phi^*\circ\iota_X$, cf. the diagram below:
\begin{center}
\includegraphics*{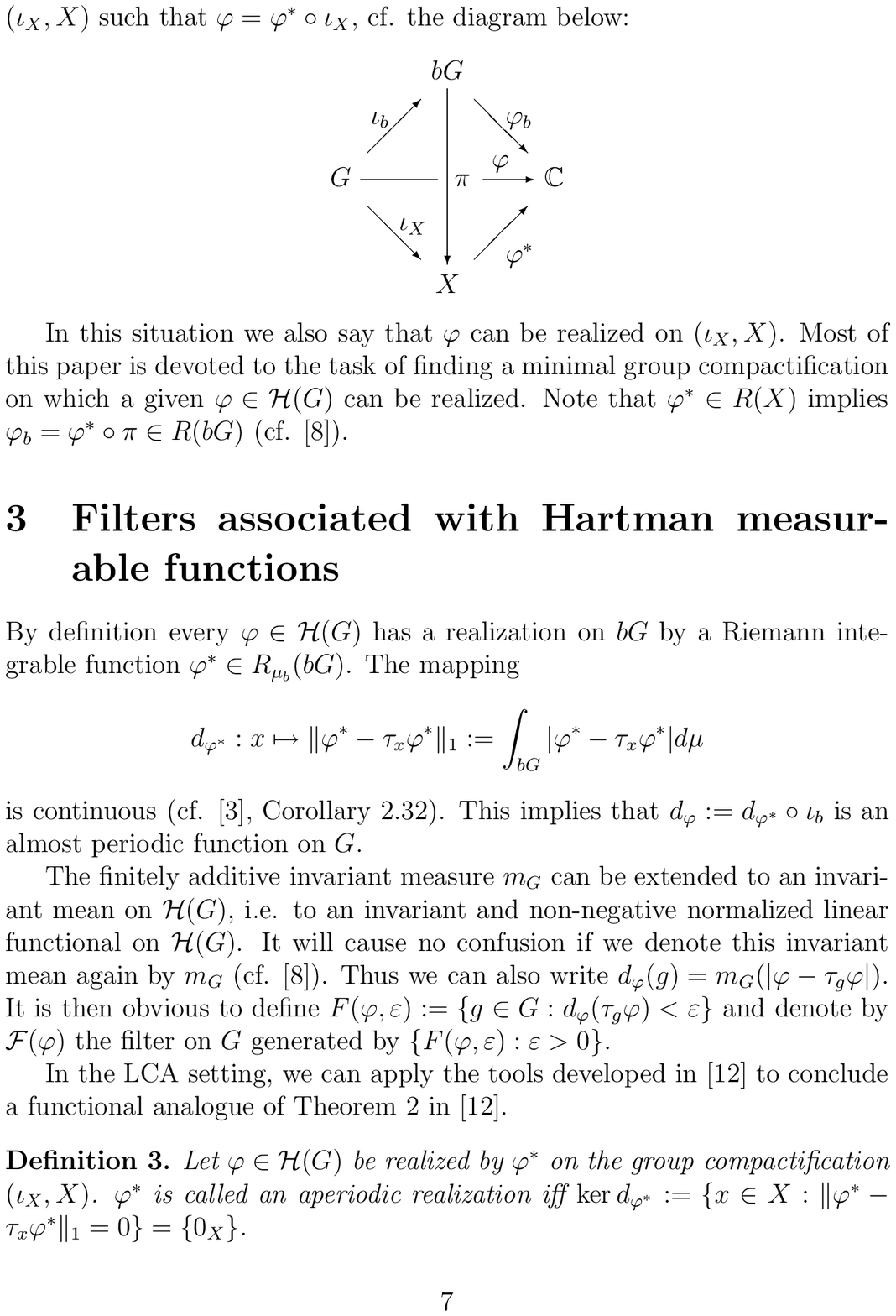}
\end{center}

In this situation we also say that $\phi$ can be realized on $(\iota_X, X)$.
Most of this paper is devoted to the task of finding a minimal group compactification
on which a given $\phi\in\H(G)$ can be realized. Note that $\phi^* \in R(X)$ implies
$\phi_b=\phi^*\circ\pi \in R(bG)$ (cf. \cite{MW}).

\section{Filters associated with Hartman measurable functions}
By definition every $\phi \in \H(G)$ has a realization on $bG$ 
by a Riemann integrable function $\phi^*\in R_{\mu_b}(bG)$.
The mapping
\[d_{\phi^*}: x\mapsto \|\phi^* -\tau_x \phi^*\|_1:=\int_{bG} |\phi^* -\tau_x \phi^*|d\mu\]
is continuous (cf. \cite{El}, Corollary 2.32). 
This implies that $d_\phi:=d_{\phi^*} \circ \iota_b$ is 
an almost periodic function on $G$. 

The finitely additive invariant measure $m_G$ can be extended to an invariant mean
on $\H(G)$, i.e. to an invariant and non-negative normalized linear functional on $\H(G)$.
It will cause no confusion if we denote this invariant mean again by $m_G$ (cf. \cite{MW}).
Thus we can also write
$d_{\phi}(g)= m_G(|\phi -\tau_g\phi|)$. It is then obvious
to define $F(\phi,\eps):=\{g\in G: d_{\phi}(\tau_g \phi)<\eps\}$
and denote by $\F(\phi)$ the filter on $G$ generated by $\{F(\phi,\eps):\eps > 0\}$.

In the LCA setting, we can apply the tools developed in \cite{Wi} to conclude a functional
analogue of Theorem 2 in \cite{Wi}.

\begin{defi}
Let $\phi\in \H(G)$ be realized by $\phi^*$ on the
group compactification $(\iota_X,X)$. $\phi^*$ is called an aperiodic
realization iff $\ker d_{\phi^*}:=\{x\in X: \|\phi^*-\tau_x \phi^*\|_1=0\}= \{0_X\}$.
\end{defi}

\begin{thr} \label{ganzleicht}
Let $\phi\in \H(G)$ be realized by $\phi^*$ on the
group compactification $(\iota_X,X)$. Then $\F(\phi)\subseteq \U_{(\iota_X,X)}$.
Furthermore  
$\F(\phi)=\U_{(\iota_X,X)}$ if $\phi^*$ is an aperiodic realization.
\end{thr}
\begin{proof}
Suppose $\phi=\phi^*\circ\iota_X$ with $\phi^*\in R_{\mu_X}(X)$ for a group
compactification $(\iota_X, X)$. For any set $A\in \F(\phi)$
there exists $\eps>0$ such that $d_{\phi}(x)<\eps$ implies $x\in A$.
Using almost periodicity
of $d_{\phi}$, i.e. continuity of $d_{\phi^*}$, we find a
neighborhood $U\in \U(X,0_X)$ such that
$d_{\phi^*}(U)\subseteq [0,\eps)$. For every 
$x\in\iota_X^{-1}(U)$ we have $d_{\phi}(x)<\eps$. Consequently 
$\iota_X^{-1}(U)\subseteq A \in \U_{(\iota_X,X)}$
and hence $\F(\phi)\subseteq \U_{(\iota_X,X)}$.

Suppose that $\phi^*\in R_{\mu_X}(X)$ is aperiodic, i.e. 
$d_{\phi^*}(x)=0$  iff $x=0_X$, the unit in $X$.
Let 
$A \in \U_{(\iota_X,X)}$ be arbitrary; w.l.o.g. we
can assume $A\supseteq \iota_X^{-1}(U)$
for an open neighborhood $U\in \U(X,0_X)$.
Due to the continuity of $d_{\phi^*}$  and compactness 
of $X$ we have $d_{\phi^*}(x)\ge\eps>0$ for $x\in X\setminus U^{\circ}$. 
This implies 
$\iota_X (\{g\in G: d_{\phi}(g)<\eps\}) \subseteq U$ and hence 
$\{g\in G: d_{\phi}(g)<\eps\}\subseteq\iota_X^{-1}(U)
\subseteq A \in \F(\phi)$.
Thus $\U_{(\iota_X,X)}\subseteq \F(\phi)$ and 
consequently $\U_{(\iota_X,X)}= \F(\phi)$.
\end{proof}

\begin{defi} Let $\phi\in \H(G)$ and let $(\iota_X, X)$ 
be a group compactification of $G$. A
function $\psi^*\in R_{\mu_X}(X)$ is called an almost
realization of $\phi$ iff
$m_G(|\phi-\psi|)=0$ for $\psi:=\psi^*\circ\iota_X$ and
$m_G$ the unique invariant mean on $\H(G)$.
\end{defi}

\begin{thr}\label{weilthr} Every $\phi\in\H(G)$ has
an aperiodic almost realization on some group compactification
$(\iota_X,X)$.
If $\phi^*: X\to \C$ is any aperiodic almost realization of $\phi$ 
then $\F(\phi)=\U_{(\iota_X,X)}$.
\end{thr}

\begin{proof}
We only have to prove that an aperiodic almost realization exists, the rest follows from
Theorem \ref{ganzleicht}.
Let $\phi^*$ be a realization of $\phi$ on $X$. The reader will easily
check that 
$H:=\ker d_{\phi^*}=\{x\in X: d_{\phi^*}(x) = 0\}$
is a closed subgroup of the compact abelian group $X$.

Weil's formula for continuous functions on quotients (Theorem 3.22 in \cite{El}) states
that there exists $\alpha>0$ such that for every $f\in C(X)$ 
\begin{eqnarray}\label{Weil}
\int_{X/H}\Big(\underbrace{\int_H f(s+t) d\mu_H(t)}_{=\b f(s)}\Big) d\mu_{X/H}(s)&=&\alpha \int_X f(u)d\mu_X(u)
\end{eqnarray} 
holds. This implies that the canonical mapping
$\b: C(X)\to C(X/H)$, $f\to \b f$ defined by $\b f(s+H)=\int_H f(s+t)d\mu_H(t)$ satisfies 
$\|\b f\|_1\le \alpha \|f\|_1$. We rescale the Haar measure on $H$ such that
$\alpha=1$.
Thus we can extend $\b$ to a continuous linear operator 
$L^1(X)\to L^1(X/H)$. Furthermore
positivity of $\b$ enables us to extend $\b$ to a mapping defined
on $R_{\mu_X}(X)$ in the following way:

According to Proposition \ref{riemannchar} $f \in R_{\mu_X}(X)$ implies that there are $g_n,h_n \in C(X)$ such
that $g_n\le f \le h_n$ and $\|h_n-g_n\|_1\to 0 $ as $n\to \infty$. Thus
every function $\tilde{f}$ on $X/H$ 
satisfying 
\begin{equation*}
f_{\bullet}:=\sup_{n\ge 0} \b g_n \le \tilde{f} \le \inf_{n\ge 0} \b h_n=:f^{\bullet}
\end{equation*} 
is in $R_{\mu_{X/H}}(X/H)$. Note that $f_{\bullet}$ and $f^{\bullet}$ are $\mu_H$-measurable
and coincide $\mu_H$-a.e.; to define $\b f$ we pick 
any function $\tilde{ f}$
satisfying $f_{\bullet}\le \tilde{f}\le f^{\bullet}$. Then Weil's formula (\ref{Weil}) will still be valid,
regardless of the particular choice of the $g_n,h_n$ and $\b f$.

Since $\phi^*$ is Riemann integrable on $X$,
there exist functions $\phi_n\in C(X)$ such that $\|\phi^*-\phi_n\|_1\to 0$.
Note that  $d_{\phi_n}\to d_{\phi^*}$ even uniformly on $X$:
\begin{equation*}
|d_{\phi_n}(s)-d_{\phi^*}(s)|  = \Big| \|\tau_s \phi_n - \phi_n\|_1 -  \|\tau_s \phi^* - \phi^*\|_1\Big| 
\le 2 \| \phi_n -\phi^*\|_1 \to 0.
\end{equation*}
Using the continuity of $\b$ as a mapping on $L^1(X)$ the same argument also shows that 
$|d_{\b \phi^*}(s+H)-d_{\b \phi_n}(s+H)|\le 2 \|\phi_n-\phi^*\|_1\to 0$ uniformly on $X/H$.
Suppose $d_{\b \phi}^*(s+H)=0$. Then
\begin{equation*}
d_{\phi^*}(s)=\lim_{n\to \infty}d_{\phi_n}(s)=\lim_{n\to \infty}d_{\b\phi_n}(s+H)=0
\end{equation*}
implies $s\in H$, i.e. $s+H=0_X+H \in X/H$. So $\b \phi^*$ is aperiodic.

We show that $\phi^*$ being a realization of $\phi$ implies that
$\b \phi^*$ is an almost realization of $\phi$.
By definition $t\in H$ iff $A_t:=\{ s \in X: \phi^*(s+t)=\phi^*(s)\}$ has $\mu_{X}$-measure $1$.
Applying Weil's formula (\ref{Weil}) to the function $f=\Eins_{A_t}\in L^1(X)$ gives
\begin{equation}\label{vorher}
\int_{X/H}\b f d\mu_{X/H}=\int_{X/H}\b \Eins_{A_t}(s+H) d\mu_{X/H}(s+H)=\int_X fd\mu_X =1.
\end{equation}
Plugging the definition of $\;\b\;$ into (\ref{vorher}) we get $\mu_{X/H}$-a.e. the identity
\begin{equation*}
\b \Eins_{A_t}(s+H)= \int_H \Eins_{A_t}(s+u)d\mu_H(u)=1.
\end{equation*}
So for every $t\in H$ and $\mu_{X/H}$-a.e. $s+H$ we know that
the set $\{u\in H: \phi^*(s+t+u)\neq\phi^*(s+u)\}$ is a $\mu_H$-null set. This
means
\[\tau_{t}(\tau_s{\phi^*}_{|H})=\tau_s{\phi^*}_{|H}\quad \mu_H\mbox{-a.e.}\]

Thus $\tau_s\phi^*$ is constant $\mu_H$-a.e. on $H$ and
for $\mu_{X/H}$ almost all $s+H$ we have
\begin{eqnarray*}
\b \phi^* (s+H)&=&\int_H \tau_s\phi^*(t) d\mu_H(t)=\int_H \phi(s)^* d\mu_H(t)=\phi^*(s).
\end{eqnarray*}

Let $\pi_H: X\to X/H$ be the quotient mapping
onto the group compactification $(\iota_{X/H},X/H)$.
Let $\psi^*:= \b\phi^*\circ \pi_H $.
Since $\b\phi^*$ is Riemann integrable on $X/H$ it is an elementary fact
that $\psi^*$ is Riemann integrable
on $X$ (cf. \cite{MW}).
Once again, Weil's formula (\ref{Weil}) together with the fact
that the Haar measure on the quotient $X/H$ is given by $\mu_{X/H}=\pi_H^{-1}(\mu_{X})$
implies
$\psi^*= \phi^*$ $\mu_{X}$-a.e. Thus the function $\psi$ defined by 
\[\psi:=\psi^*\circ \iota_X = \b \phi^*\circ\iota_{X/H}\] satisfies 
$m_G(|\phi-\psi|)=\|\phi^*-\psi^*\|_1=0$ for the unique invariant mean $m_G$.
Thus $\psi^*$ is the required almost realization of $\phi$.
\end{proof}

\begin{cor}
\label{apreal}
Every $\phi\in\A(G)$ has an aperiodic realization 
on some group compactification $(\iota_X,X)$.
\end{cor}
\begin{proof}
We use the notation from Theorem \ref{weilthr}.
If $\phi$ is almost periodic then $\phi^*$ is continuous.
Consequently $\b \phi^*$ and
$\psi^*:=\b\phi^*\circ\pi$ are also continuous. Since these functions
coincide $\mu_X$-a.e. 
they coincide everywhere on $X$. This implies that
$\phi^*$ is constant on $H$-cosets and
$\b \phi^* (s+H)=\phi^*(s)$ for all $s+H \in X/H$. So $\phi^*$ is truly
a realization of $\phi$, not only an almost realization.
\end{proof} 

This Corollary is a special case of F\o lner's
"Main Theorem for Almost Periodic Functions", for a detailed treatment
cf. \cite{Fo}.

\emph{Remark:} Note that for any given realization of a Hartman measurable
function $\phi\in\H(G)$ on a group compactification $(\iota_X,X)$
we can always assume that there exists an aperiodic almost realization of $\phi$ on a
group compactification $(\tilde{\iota_X},\tilde{X})$ with 
$(\tilde{\iota_X},\tilde{X})\le (\iota_X,X)$.
Since in \cite{MW} it is shown
that every Hartman measurable function on an LCA group
with separable dual has a realization on  a
metrizable group compactification,
every Hartman measurable function 
on such a group has
an aperiodic almost realization on a metrizable group compactification.

\begin{lem}\label{compsub} Let $G$ be an LCA group and let $(\iota_X,X)$ 
be a group
compactification. Then there exists a unique subgroup $\Gamma\le \hat{G}$
such that $(\iota_{\Gamma},C_{\Gamma})$ and $(\iota_X,X)$ are equivalent.
Furthermore $(\iota_X,X)$ is the supremum of 
all group compactifications $(\iota_{\gamma},C_{\gamma})$
such that $(\iota_{\gamma},C_{\gamma})\le (\iota_X,X)$ (writing  
in short $(\iota_{\gamma},C_{\gamma})$ for
$(\iota_{\langle\gamma\rangle},C_{\langle\gamma\rangle})$).

The mapping $(\iota_X,X)\mapsto C_{\Gamma}$ is a bijection between
equivalence classes of group compactifications of $G$
and subgroups of $\hat{G}$.
\end{lem} \begin{proof} See Theorem 26.13 in \cite{HR}. \end{proof}

\begin{cor}\label{unique} Let $\phi\in\H(G)$. Any two group compactifications 
$(\iota_{X_1},X_1)$ and $(\iota_{X_2},X_2)$ on which
$\phi$ has an aperiodic almost realization are equivalent. \end{cor}
\begin{proof} By Theorem \ref{ganzleicht} we have $\U_{(\iota_{X_1},X_1)}=\F(\phi)
=\U_{(\iota_{X_2},X_2)}$.
A straight forward generalization of Theorem 1 in \cite{Wi} implies
that the mapping $$\Phi: \hat{G}\ge \Gamma \mapsto (\iota_{{\Gamma}},C_{\Gamma})$$ coincides
with the composition of the mappings
\begin{eqnarray*}
\Sigma: (\iota_b,bG)\ge(\iota_X,X)  &\mapsto& \U_{(\iota_X,X)},\\
\Psi: \mathfrak{P}(G)\supset \F   &\mapsto& 
\{\chi \in \hat{G}: \F\!\!-\!\!\lim_{g\in G}\chi(g) = 0\}.
\end{eqnarray*}
Since Lemma \ref{compsub} states that $\Phi=\Psi\circ\Sigma$ is invertible,
$\Sigma$ must be one-one. In particular 
$\U_{(\iota_{X_1},X_1)}=\U_{(\iota_{X_2},X_2)}$ implies that
$(\iota_{X_1},X_1)$ and $(\iota_{X_2},X_2)$ are 
equivalent group compactifications.
\end{proof}

{\em For the rest of this section assume that $G$ is an 
LCA group with separable dual.}

\begin{cor}
Every filter $\F(\phi)$ with $\phi\in\H(G)$ coincides with a filter 
$\U_{(\iota_X,X)}$ for a metrizable
group compactification $(\iota_X,X)$. If $\phi^*$ is 
an arbitrary realization of $\phi$, say on the Bohr compactification
$bG$, we can take $X\cong bG/\ker d_{\phi^*}$.
\end{cor}

\begin{cor}
Hartman measurable functions 
induce exactly the filters coming from metrizable group compactifications.
\end{cor}
\begin{proof} In Theorem 3 in \cite{Wi} 
for every metrizable group compactification $(\iota_X,X)$ of the 
integers $G=\Z$, an aperiodic Hartman periodic function of the form
$f=\Eins_A$ is constructed. The same construction can be done in an arbitrary
LCA group $G$ as long as the dual $\hat{G}$ contains a countable and dense
subset.
This shows that any $\U_{(\iota_X,X)}$ with metrizable $X$ can be obtained
already by Hartman measurable \emph{sets},
i.e. by a filter $\F(\phi)$ with $\phi=\Eins_A$. 
Since 
any Hartman measurable function on $G$ can be realized on a metrizable group compactification
(cf. \cite{MW}).
Thus Theorem \ref{ganzleicht} implies that
no filter $\F(\phi)$ can coincide with $\U_{(\iota_X,X)}$ for a 
non metrizable group compactification $(\iota,C)$.
\end{proof}

\section{Subgroups associated with Hartman measurable functions}

For Hartman measurable $\phi$
let us denote by $\Gamma(\phi)$ the (countable) subgroup of $\hat{G}$ generated by 
the set 
\[\mbox{spec}\;\phi:=\{\chi\ \in \hat{G}: m_G(\phi\cdot\overline{\chi})\neq 0\}\]
of all characters with
non vanishing Fourier coefficients. We will prove that $\Gamma=\Gamma(\phi)$
determines
a group compactification $(\iota_{\Gamma}, C_{\Gamma})$ such that $\phi$ can
be realized aperiodically on $C_{\Gamma}$. First we deal with almost periodic functions:

\begin{pro}\label{fact1} Let $\phi\in\A(G)$ and 
$(\iota_X, X)$  a group compactification
such that
every character $\chi\in \Gamma(\phi)$ has a representation $\chi = \eta \circ \iota_X$ with a
continuous character $\eta \in \hat{X}$.
Then every function $f\in \overline{\mbox{span}}\,\Gamma(\phi)\subseteq AP(G)$
has a realization on $(\iota_X, X)$. \end{pro}
\begin{proof}
This is essentially a reformulation of Theorem 5.7 in \cite{B}.
In fact the  Stone-Weierstrass Theorem implies
that $\overline{\mbox{span}}\,\Gamma(\phi)= \iota_{\Gamma}^*\;C(X)$. Furthermore 
$\phi\in\overline{\mbox{span}}\,{\Gamma}(\phi)$, i.e. $\phi$ can be realized by some
continuous $\phi^*: X \to \C$. \end{proof}

\begin{pro}\label{characterreal1} Let $\phi\in\A(G)$ and 
$(\iota_{\Gamma}, C_{\Gamma})$
the group compactification of $G$ induced
by the subgroup $\Gamma=\Gamma(\phi)\le \hat{G}$.
Then for every continuous character $\psi\in \Gamma(\phi)$ 
there exists a continuous $\psi^*: C_{\Gamma} \to \C$ such that
$\psi=\psi^*\circ\iota_{\Gamma}$. \end{pro}
\begin{proof}

Given the group compactification $(\iota_{\Gamma},C_{\Gamma})$, then
the compact group $C_{\Gamma}$ is 
by definition topologically isomorphic to
$\overline{\{(\chi(g))_{\chi\in \Gamma}: g\in G\}} \le 
\T^{\Gamma}$.

The restriction of each projection 
\[\pi_{\chi_0}: C_{\Gamma}\le\T^{\Gamma} \to \T,\quad 
(x_{\chi})_{\chi \in \Gamma}\mapsto x_{\chi_0} \]
is a bounded character of $C_{\Gamma}$ for each $\chi_0 \in \Gamma(\phi)$.
I.e. $\pi_{\chi_0}$ is an element of $\widehat{C_{\Gamma}}$. 
Thus $\chi_0= \pi_{\chi_0} \circ \iota_{\Gamma}$ for each 
$\chi_0\in \Gamma(\phi)$ and
we may apply Proposition \ref{fact1} to obtain the assertion.
\end{proof}

\begin{pro}\label{characterreal2}
Let $\phi\in \A(G)$ and let 
$(\iota_X,X)$ be a
group compactification of $G$ such
that $\phi$ can be realized by a continuous function $\phi^*: X \to \C$. 
Then each continuous character $\chi\in \Gamma(\phi)$ has a representation
$\chi = \eta \circ \iota_{\Gamma}$ with $\eta \in \hat{X}$. \end{pro}

\begin{proof} Obviously it is enough to prove the assertion for 
a generating subset of $\Gamma(\phi)$. Let
$\chi\in \hat{G}$ be such
that $m_G(\phi\cdot\overline{\chi})\neq 0$. Define a linear functional 
$m_{\chi}: C(X) \to \C$
via $\psi\mapsto m_{\chi}(\psi)=m_G(( \psi\circ\iota_{\Gamma}) \cdot \overline{\chi})$.
It is routine to check
that $m_{\chi}$ is bounded and $\|m_{\chi}\|=1$. Since $X$ is compact the 
complex-valued mapping $\tilde{\eta}: 
X \mapsto m_{\chi}(\tau_x \phi^*)$ 
is continuous on $X$ (the mapping $x\mapsto \tau_x \phi^*$ is continuous).
For $g\in G$ we compute
\begin{eqnarray*}
\tilde{\eta} \circ \iota_X (g) &=& m_G((\tau_{\iota_X(g)} \phi^*\circ \iota_X) \cdot\overline{\chi}) 
                                = m_G(\tau_g (\phi^*\circ \iota_X)\cdot\overline{\chi})\\
                               &=& m_G( (\phi^*\circ \iota_X)\cdot \tau_{-g}\overline{\chi})
							    = m_G( (\phi^*\circ \iota_X) \cdot\chi(g)\overline{\chi})\\
							   &=&\chi(g) m_{\chi}(\phi^*)=  \chi(g) \tilde{\eta} (0).
\end{eqnarray*}
Since $\tilde{\eta}(0)=m_{\chi}(\phi^*)=m_G(\phi \cdot \overline{\chi})\neq 0$ we can define
$\eta := \tilde{\eta}(0)^{-1}\tilde{\eta}$. The mapping
$\eta: X \to \T$ is continuous
and satisfies the functional equation 
\[
\eta(\iota_X(g)+\iota_X(h))=
\tilde{\eta}(0)^{-1}\tilde{\eta}(\iota_X(g)+\iota_X(h))=
\chi(g)\chi(h)=\eta(\iota_X(g))\eta(\iota_X(h))\]
on the dense set $\iota_X(G)$.
Hence $\eta$ is a bounded character on $X$
and $\eta\circ\iota_X = \chi$. 
\end{proof}

\begin{cor}\label{characterreal2a}
Let $\phi\in \H(G)$ be realized by $\phi^*$ on the
group compactification $(\iota_X,X)$.
Then each $\chi\in \Gamma(\phi)$ has a representation
$\chi = \eta \circ \iota_X$ with $\eta \in \hat{X}$. \end{cor}
\begin{proof} For every $\chi\in\hat{G}$ with 
$m_G(\phi\cdot\overline{\chi})=\alpha\neq 0$
we can pick a continuous function $\psi^*:X\to \C$
such that $\|\psi^*-\phi\|_1 < |\alpha|/2$. Then
$\psi:=\psi^*\circ\iota_X$ satisfies 
\[|m_G(\phi\cdot\overline{\chi})-m_G(\psi\cdot\overline{\chi})|\le m_G(|\phi-\psi|)\le 
\|\psi^*-\phi^*\|_1 < |\alpha|/2.\] In particular $m_G(\psi\cdot\overline{\chi})\neq 0$.
Applying Proposition \ref{characterreal2} to the function $\psi\in\A(G)$ yields that
the character $\chi$ can be realized on $X$.\end{proof}

Thus for almost periodic functions $\phi$ the subgroup $\Gamma(\phi)$ contains all the 
relevant information to reconstruct $\phi$ from its Fourier-data in a minimal way.
It is not obvious how to obtain similar results for 
Hartman measurable functions that are not
almost periodic. The following example illustrates how a straight forward
approach may fail.

\begin{exa}{\em
Let $\phi_n(k):=\prod_{j=1}^n \cos^2\left(2\pi \frac{k}{3^j}\right)$
on $G=\Z$. Each $\phi_n$ is a finite product of periodic (and hence almost periodic)
functions. Since $\A(\Z)$ is an algebra, $\phi$ is almost periodic.
In \cite{MW} it is shown that
$\phi(k):=\lim_{n\to \infty} \phi_n(k)$ exists and defines
a non negative Hartman measurable function with $m_{\Z}(\phi)=0$.
Since $\Gamma(\phi_n)\cong\Z/3^n \Z$ we have (using obvious notation):
\[\lim_{n\to\infty} \Gamma(\phi_n)=
\bigcup_{n=1}^{\infty} \Gamma(\phi_n)\cong
\Z_3^{\infty},\] the Pr\"ufer 3-group (i.e. the subgroup of all
complex $3^n$-th roots of unity for $n\in \N$), 
but
\[\Gamma(\lim_{n\to\infty} \phi_n)=\{0\}.\]}\end{exa}

\begin{pro}\label{fejer}
Let $\{K_n\}_{n=1}^{\infty}$ denote the family of
Fej\'er kernels on $\T^k$
\[ K_n(\exp(it_1),\ldots, \exp(it_k))=\frac{1}{k}
\prod_{j=1}^{k} \left(\frac{\sin( \frac{1}{2}nt_j)}{\sin( \frac{1}{2}t_j)}\right)^2.\]
The linear convolution operators on $L^1(\T^k)$ defined
by \[\sigma_n: \phi\mapsto K_n*\phi\] are
non negative, their norm is uniformly bounded by 
$\|\sigma_n\|=1$ and $\sigma_n \phi(x) \to \phi(x)$ a.e. for every $\phi\in L^1(\T^k)$.
Furthermore $\sigma_n \phi \in \mbox{span } \Gamma(\phi)$ for every $n\in \N$.
\end{pro} \begin{proof} This is a reformulation
of the results in section 44.51 in \cite{HR}.\end{proof}

Let $f$ be Riemann integrable on $X=\T^k$, w.l.o.g.
real-valued, and $\phi_{i}, \psi_i\in C(X)$ such that
$\phi_{i}\ge f\ge \psi_i$ and $\|\phi_{i}-\psi_i\|_1<\eps_i$ for a sequence 
$\{\eps_i\}_{i=1}^{\infty}$ of positive real numbers, tending monotonically to $0$.
We know that $\sigma_n f(x)\to f(x)$ for a.e. $x\in X$. Thus we have
\begin{equation*}\label{rapf}
\phi_n^*:=\sigma_n \phi_n \ge \sigma_n f \ge \sigma_n \psi_n=\psi_n^*
\end{equation*}
and
\begin{equation*}
\|\phi_n^*-\psi_n^*\|_1\le \|\sigma_n (\phi_n^*-\psi_n^*)\|_1
                       \le  \|\sigma_n\|\, \|\phi_n-\psi_n\|_1 \le \eps_n.
\end{equation*}
Let $\phi^*:=\inf_{n\in\N}\phi_n$ and 
$\psi^*:=\sup_{n\in\N}\psi_n$. If we assume w.l.o.g. $\psi_n$ to increase and $\phi_n$ to decrease as $n\to\infty$, 
the same will hold for $\psi_n^*$ and $\phi_n^*$. This implies that in the inequality 
\[\phi^*(x)=\lim_{n\to\infty} \phi_n^*(x) \ge \limsup_{n\to\infty} \sigma_n f
\ge \liminf_{n\to\infty} \sigma_n f \ge
\lim_{n\to\infty}\psi_n^*(x)=\psi^*(x)\]
actually equality holds $\mu_X$-a.e. on $X$.  Thus we can apply
Proposition \ref{riemannchar} and conclude that any function $f^*$
with $\phi^*\ge f^*\ge\psi^*$ is Riemann integrable (and coincides
$\mu_X$-a.e. with $f$). 
In particular $f^{\bullet}:=\limsup_{n\to \infty} \sigma_n f$ and
$f_{\bullet}:=\liminf_{n\to \infty} \sigma_n f$
are (lower resp. upper semicontinuous) Riemann integrable 
functions that coincide $\mu_X$-a.e. with $f$.

Let us call 
a group compactification $(\iota_X,X)$ finite dimensional iff
$X$ is topologically isomorphic
to a closed subgroup of $\T^n$ for some $n\in \N$.
Note that if $(\iota_X,X)$ is finite dimensional, then every
group compactification covered by $(\iota_X,X)$ is finite dimensional
as well.
 A Hartman measurable function $\phi\in\H(G)$ can be realized
finite dimensionally iff there exists a realization of $\phi$ on some
finite dimensional group compactification.

\begin{pro}
For a compact LCA group $C$ the following assertions are equivalent:
\begin{enumerate}
\item $C$ is finite dimensional,
\item $\hat{C}$ is finitely generated,
\item $C$ is topological isomorphic to $\T^k\times F$ for $k\in \N$
and a finite group $F$ of the
form \[F=\prod_{i=1}^N (\Z/n_i\Z)^{p_i},\quad p_i \mbox{ prime}.\]
\end{enumerate}

\end{pro}
\begin{proof} Folklore.
\end{proof}

\begin{pro}
Let $\phi\in\H(G)$. If $\phi$ can be realized finite dimensionally, then there
is an almost realization of $\phi$ on the (finite dimensional) compactification
induced by $\Gamma:=\Gamma(\phi)$.
\end{pro}

\begin{proof}

Let $\phi$ be realized finite dimensionally on some group compactification
$(\iota_X,X)$. Since there exists a group compactification 
covered by $(\iota_X,X)$, on which $\phi$ can be almost realized aperiodically (cf. 
Theorem \ref{weilthr}), we can assume w.l.o.g. that $\phi$ can be almost realized aperiodically
already on $(\iota_X,X)$. We have to show that $(\iota_X,X)$ and 
$(\iota_{\Gamma},C_{\Gamma})$ are equivalent.

Let $\psi^*$ be an aperiodic realization of $\phi$ on 
$C_{\Gamma}\cong \T^k\times F$ with $k\in\N$ and $F$ finite.
Let us denote the elements of $\T^k\times F$ by tuples $(\vec{\alpha},x)$.
For every fixed $\vec{\alpha}\in \T^k$ define 
a mapping $\psi_{\vec{\alpha}}: F\to \R$ via
\[\psi^*_{\vec{\alpha}}(x):=\psi^*(\vec{\alpha},x).\]
For each $\overline{\chi} \in \hat{F}$, the dual of the 
finite group $F$, define the $F$-Fourier coefficient of $\psi^*_{\vec{\alpha}}$ as
\[c_{\overline{\chi}}(\vec{\alpha}):= \int_F \psi^*_{\vec{\alpha}}(x)\overline{\chi}(x) dx= \frac{1}{\# F} 
\sum_{x\in F} \psi^* (\vec{\alpha};x) \overline{\chi}(x)\in\C.\]
We want to show that $c_{\overline{\chi}}: \T^k\to \C$ is a Riemann integrable function: 
The mapping $\gamma_x: \T^k\to X$ defined via
$\vec{\alpha}\mapsto(\vec{\alpha};x)$ is continuous and measure-preserving
for every $x\in F$. $\psi^*$ is by definition Riemann integrable. Thus the mapping 
$\psi^*\circ\gamma_x: \T^k\to \C $ is Riemann integrable for each $x\in F$.
Note that
\[c_{\overline{\chi}}(\vec{\alpha})=
\sum_{x\in F} (\psi^*\circ\gamma_x) (\vec{\alpha}) \overline{\chi}(x).\]
Hence, for each fixed character $\chi \in \hat{F}$,
the mapping $c_{\overline{\chi}}: \T^k \to \C$ defined 
via $\vec{\alpha}\mapsto \sum_{x\in F} (\psi^*\circ\gamma_x) (\vec{\alpha}) \overline{\chi}(x)$
is Riemann integrable on $\T^k$.

Thus Proposition \ref{fejer} implies
$\sigma_n c_{\chi}(\vec{\alpha})\to c_{\chi}(\vec{\alpha}) $ a.e. on $\T^k$.
Taking into account 
that the Haar measure on $F$ is the normalized counting measure, we get
\begin{equation}\label{syn} \psi^*_n(\vec{\alpha};x):=
\sum_{\overline{\chi}\in \hat{F}} \left(\sigma_n c_{\overline{\chi}}(\vec{\alpha})\right)
\overline{\chi}(x) \to \sum_{\overline{\chi}\in \hat{F}} c_{\overline{\chi}}(\vec{\alpha})\overline{\chi}(x) = 
\psi^*_{\vec{\alpha}}(x)= \psi^*(\vec{\alpha};x)\end{equation}
for almost every $\vec{\alpha}\in \T^k$ and every $x\in F$, as $n\to \infty$.
Since Haar measure $\mu_C$ on $C$ is the product measure
of the Haar measures on the groups $\T^k$ and $F$, the relation (\ref{syn}) 
holds $\mu_C$-a.e. on $C$.
We conclude that
any function majorizing $\liminf_{n\to \infty} \psi^*_n$ and minorizing $\limsup_{n\to \infty} \psi^*_n$
is an almost realization of $\phi$.
Note that according to the properties of the Fej\'er kernels on $\T^k$ 
(see 44.51 in \cite{HR}) for each character
$(\eta\times\chi) (\vec{\alpha};x):=\eta(\vec{\alpha})\chi(x)$, 
$\eta\in \hat{\T^k}$ and $\chi\in\hat{F}$, 
there exists an $n_0\in N$ such that for $n\ge n_0$
in the Fourier expansion of $\psi^*_n$ 
the Fourier coefficient (computed in $C$) associated with the character
does not vanish iff the $\T^k$-Fourier coefficient of $c_{\chi}$
\[c_{\eta}(c_{\chi})=\int_{\T^k} c_{\chi }(\vec{\alpha})\overline{\eta}(\vec{\alpha})d\vec{\alpha}\]
does not vanish.
A simple computation shows that
the Fourier coefficients of $\psi^*$ are given by 
\begin{eqnarray*}
c_{\eta \times \chi}(\psi^*)&=&
\int_{\T^k}\int_F \psi^*(\vec{\alpha},x)\overline{\eta}(\vec{\alpha})\overline{\chi}(x) d\vec{\alpha} dx \\
&=& \int_{\T^k} c_{\chi}(\vec{\alpha}) \overline{\eta}(\vec{\alpha}) d\vec{\alpha}
= c_{\eta}(c_{\chi}) \end{eqnarray*}

So the character $\eta\times\chi$ contributes to the Fourier
expansion of $\psi^*$ if and only if
$c_{\eta \times \chi}(\psi)\neq 0$. Thus
$\psi^*_n\in \mbox{span } \Gamma(\phi)$ for every $n\in\N$,
implying that there exist almost realizations of $\phi$ on the
group compactification
induced by $\Gamma(\phi)$, e.g. $\liminf_{n\to \infty} \psi^*_n$
or $\limsup_{n\to \infty} \psi^*_n$.
\end{proof}

Combining this result with the results of the previous section
we obtain
\begin{thr}\label{mincorap}
Let $\phi\in\H(G)$ and $\Gamma=\Gamma(\phi)\le \hat{G}$. The following assertions hold:
\begin{enumerate}
\item $(\iota_{\Gamma},C_{\Gamma})\le (\iota_X,X)$ for every
compactification $(\iota_X,X)$ on which $\phi$ can be realized.
In particular $\F(\phi)\subseteq \U_{(\iota_{\Gamma},C_{\Gamma})}$.
\item Assume that $\phi\in \A(G)$ or that $\phi$ can be realized finite
dimensionally. Then $\phi$ can be realized aperiodically
on $C_{\Gamma}$. In particular $\F(\phi)=\U_{(\iota_{\Gamma},C_{\Gamma})}$.
\end{enumerate}
\end{thr}

We strongly conjecture that the second assertion in Theorem \ref{mincorap} 
holds for \emph{any}
Hartman measurable function, at least on LCA groups $G$ with 
separable dual $\hat{G}$. A proof of this might
utilize more general summation methods
(in the flavour of Theorems 44.43 and 44.47 in \cite{HR}) than the Fej\'er summation
presented here.

In \cite{Wi} it is shown that for any Hartman measurable set $M\subseteq G=\Z$ and the
induced filter $\F=\F(M)$ there is an aperiodic
realization of $\phi_M=\Eins_M$ on the compactification determined by the subgroup
Sub$(M)=\{\alpha: \F\!\!-\!\!\lim_{n\in\Z} \lfloor n\alpha\rfloor=0\}$ or, equivalently,
Sub$(M)=\{\alpha: \F\!\!-\!\!\lim_{n\in\Z} e^{2\pi i n\alpha}=1\}$.

Together with
Theorem \ref{mincorap} this implies that for Hartman sets $M$
with finite dimensional realization 
both the group compactifications of $\Z$ induced by the
subgroups $\Gamma(\phi_M)$ and Sub$(M)$ admit aperiodic realizations of $\phi_M$. 
Hence uniqueness of the minimal compactification with 
aperiodic realization
(Corollary \ref{unique})
implies that in this special case $\Gamma(\phi_M)=\mbox{Sub}(M)$.
In the general situation we can prove up to now far only the following

\begin{pro}
For a Hartman measurable function $\phi\in\H(G)$ let
$\F=\F(\phi)$, $\Gamma=\Gamma(\phi)$ and $\mbox{Sub}(\phi)=\{\chi\in \hat{G}:
\F\!\!-\!\!\lim_{g\in G} \chi(g) = 1_{\C}\}$.
Then $\Gamma(\phi)\le \mbox{Sub}(\phi)$.
\end{pro}
\begin{proof}
Suppose $\chi \in \Gamma(\phi)$. To prove $\F\!\!-\!\!\lim_{g\in G} \chi(g) = 1_{\C}$
(unit element of the multiplicative group of complex numbers)
we have to show that for every $\eps>0$ the set
$\{g\in G: |1-\chi(g)|<\eps\}$ belongs to the filter $\F(\phi)$,
i.e. that there exists $\delta=\delta(\eps)>0$ such that
\begin{equation}\label{filterlimes}
\{g\in G: m_G(|\tau_g\phi-\phi|)<\delta\}\subseteq \{g\in G: |1-\chi(g)|<\eps\}\in\F(\phi).
\end{equation}

Using the fact that $m_G$ is an invariant mean and that $\chi$ is a homomorphism,
we have
\[
\chi(g)\,m_G(\phi\cdot \overline{\chi})= m_G(\tau_g\phi\cdot \overline{\chi})=
m_G( (\tau_g\phi-\phi)\cdot\overline{\chi})+ m_G(\phi\cdot \overline{\chi}).
\]
Using $\|\chi\|_{\infty}=1$ this implies
\begin{equation*}\label{fcesti}
|1-\chi(g)|\cdot |m_G(\phi\cdot \overline{\chi})| = |m_G( (\tau_g\phi-\phi)\cdot\overline{\chi})|
\le m_G(|\tau_g\phi-\phi|).
\end{equation*}
Since $m_G(\phi\cdot \overline{\chi})\neq 0$ we can define
$\delta:= \eps\cdot \frac{m_G(|\tau_g\phi-\phi|)}{|m_G(\phi\cdot \overline{\chi})|}>0$.
With this choice of $\delta$ indeed
$m_G(|\tau_g\phi-\phi|)<\delta$ implies
$|1-\chi(g)|<\eps$, i.e. the inclusion (\ref{filterlimes}) holds.
\end{proof} 

\section{Summary}
The content of the present paper essentially 
deals with the  definition and properties of the objects occurring
in the diagram below. Abusing the terminus technicus of \emph{commutative
diagrams} in a kind of sloppy way, the theorems of this paper
circle around the question under which assumptions this diagram is commutative:
\begin{center}
\includegraphics*{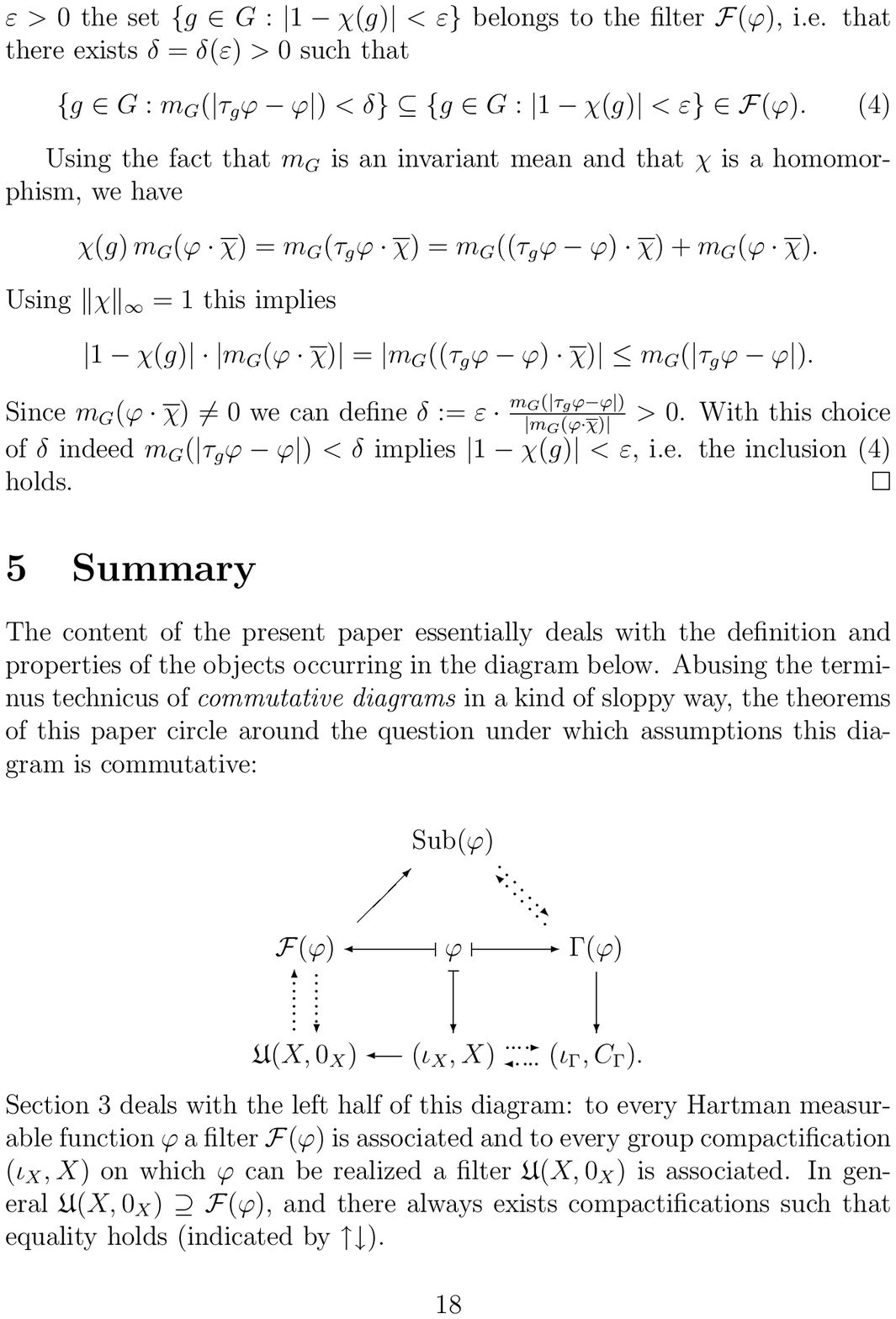}
\end{center}
Section 3 deals with the left half of this diagram: to every Hartman measurable function
$\phi$ a filter $\F(\phi)$ is associated and to every group compactification $(\iota_X,X)$ on
which $\phi$ can be realized a filter $\U(X,0_X)$ is associated. In general $\U(X,0_X)\supseteq
\F(\phi)$, and there always exists compactifications such that equality holds (indicated by
$\uparrow\downarrow$).

Section 4 deals with the right half of the diagram: to every Hartman measurable function
$\phi$ a subgroup $\Gamma(\phi)$ of the dual is associated, which in turn induces a group
compactification $(\iota_{\Gamma},C_{\Gamma})$. In general 
$(\iota_{\Gamma},C_{\Gamma})\le (\iota_X,X)$ for every  group compactification $(\iota_X,X)$ on
which $\phi$ can be realized. If $\phi$ is either almost periodic or can be realized finite
dimensionally then $(\iota_A,X_A)$ is itself a group compactification on
which $\phi$ can be realized (indicated by
{\tiny${\rightarrow\atop\leftarrow}$}) and the filter $\U(X_A,0_{X_A})$ associated with this particular
compactification coincides with $\F(\phi)$.
The filter $\F(\phi)$ in turn defines a subgroup Sub($\phi$) of the dual $\hat{G}$. While it
can be shown that in general $\Gamma(\phi)\le\mbox{Sub}(\phi) $ it is and open problem
whether this inclusion can be reversed.

\vspace{1em}\begin{center}
Author's address:\\
Gabriel Maresch\\
Technical University Vienna\\
Institute of Discrete Mathematics and Geometry\\
Wiedener Hauptstra\ss e 8-10\\
1040 Vienna, Austria\\
email: {\ttfamily gabriel.maresch@tuwien.ac.at}\\
web: http://www.dmg.tuwien.ac.at/maresch\end{center}

\begin{thebibliography}{AA}
\bibitem{B} R. Burckel, \emph{Weakly Almost Periodic Functions on Semigroups}, Gordon
and Breach, New York, 1970.
\bibitem{BS} M. Beiglb\"ock, C. Steineder, R. Winkler
\emph{Sequences and filters of characters characterizing subgroups
of compact abelian groups}, to appear in Top. Appl.
\bibitem{El} J. Elstrodt, \emph{Ma\ss- und Integrationstheorie}, Springer-Verlag
 Berlin Heidelberg New York, 1999.
\bibitem{Fo} E. F\o lner, \emph{A Proof of the Main Theorem for Almost Periodic
Functions in an Abelian Group}, Ann. of Math., 50/5 (1949), 559-569.
\bibitem{FP} S. Frisch, M. Pa\v{s}teka, R. Tichy, R. Winkler, \emph{Finitely additive measures
on groups and rings}, Rend. Circ. Mat. Palermo, Series II, 48 (1999), 323-340.
\bibitem{H} S. Hartman, \emph{Remarks on equidistribution on non-compact groups},
            Compositio Math., 16 (1964), 66-71.
\bibitem{HR} E. Hewitt and K. Ross, \emph{Abstract Harmonic Analysis I,II}, 
Springer-Verlag Berlin-Heidelberg-New York, 1963.
\bibitem{MW} G. Maresch and R. Winkler, \emph{Hartman measurable functions and 
related function spaces}, E-print, 2005, available at \emph{www.dmg.tuwien.ac.at/maresch}
\bibitem{SchW} J. Schmeling, E. Szab\`o, R. Winkler, \emph{Hartman and Beatty bisequences},
Algebraic Number Theory and Diophantine analysis, 405-421, Walter de Gruyter, Berlin New York, 2000.
\bibitem{StW} C. Steineder and R. Winkler, \emph{Complexity of Hartman sequences},
to appear in Journal de Th\'{e}orie des Nombres de Bordeaux.
\bibitem{Ta} M. Talagrand, \emph{Closed convex hull of measurable functions,
Riemann measurable functions and measurability of trans\-lations},
 Ann. Inst. Fourier (Grenoble), 32/1 (1982), 39-69.
\bibitem{Wi} R. Winkler, \emph{Ergodic Group Rotations, Hartman Sets and Kronecker
Sequences}, Monatsh. Math., 135 (2002), 333-343.
\end{thebibliography}
\end{document}